\documentclass[12pt]{amsart}

\usepackage{amsmath,amssymb,amsthm}
\usepackage{graphicx}
\usepackage{enumerate}
\usepackage{color}
\usepackage[german, english]{babel}
\usepackage[all]{xy}

\newtheorem{thm}{Theorem}[]
\newtheorem*{thm*}{Theorem}
\newtheorem{lem}[thm]{Lemma}

\newtheorem{cor}[thm]{Corollary}

\newtheorem{ex}[thm]{Example}

\newcommand{\param}{{\mathchoice{\mkern1mu\mbox{\raise2.2pt\hbox{$
\centerdot$}}
\mkern1mu}{\mkern1mu\mbox{\raise2.2pt\hbox{$\centerdot$}}\mkern1mu}{
\mkern1.5mu\centerdot\mkern1.5mu}{\mkern1.5mu\centerdot\mkern1.5mu}}}

\renewcommand{\setminus}{{\smallsetminus}}

\renewcommand \color [2][]{}

\begin{document}

\title       {Geometry of the Homology Curve Complex}
\author   {Ingrid Irmer}
\email      {irmer@cs.uni-bonn.de}
\maketitle
  
\begin{abstract}
Suppose $S$ is a closed, oriented surface of genus at least two. This paper investigates the geometry of the homology multicurve complex, $\mathcal{HC}(S,\alpha)$, of $S$; a complex closely related to complexes studied by Bestvina-Bux-Margalit and Hatcher. A path in $\mathcal{HC}(S,\alpha)$ corresponds to a homotopy class of immersed surfaces in $S\times I$. This observation is used to devise a simple algorithm for constructing quasi-geodesics connecting any two vertices in $\mathcal{HC}(S,\alpha)$, and for constructing minimal genus surfaces in $S\times I$. It is proven that for $g \geq 3$ the best possible bound on the distance between two vertices in $\mathcal{HC}(S, \alpha)$ depends linearly on their intersection number, in contrast to the logarithmic bound obtained in the complex of curves. For $g \geq 4$ it is shown that $\mathcal{HC}(S, \alpha)$ is not $\delta$-hyperbolic. 
\end{abstract}

\section{Introduction} 
Suppose $S$ is a closed oriented surface. $S$ is not required to be connected but every component is assumed to have genus $g\geq 2$.  

Let $\alpha$ be a nontrivial element of $H_{1}(S,\mathbb{Z})$. The \textit{homology curve complex}, $\mathcal{HC}(S,\alpha)$, is a simplicial complex whose vertex set is the set of all homotopy classes of oriented multicurves in $S$ in the homology class $\alpha$. A set of vertices $m_{1},\ldots, m_k$ spans a simplex if there is a set of pairwise disjoint representatives of the homotopy classes. 

The \textit{distance}, $d_{\mathcal{H}}(v_{1},v_{2})$, between two vertices $v_1$ and $v_2$ is defined to be the distance in the path metric of the one-skeleton, where all edges have length one.


The \textit{Torelli group} is the subgroup of the mapping class group that acts trivially on homology. $\mathcal{HC}(S, \alpha)$ is closely related to a complex defined in \cite{BBM} that was used for calculating cohomological properties of the Torelli group.

Metric properties of curve complexes have been used for example for studying mapping class groups and the structure of 3-manifolds, for example \cite{endinglamination}, \cite{MasurandMinskyII} and \cite{HempelII}. The aim of this paper is to study some basic geometric properties of $\mathcal{HC}(S,\alpha)$.

In \cite{MasurandMinskyI} and \cite{Bowditch2} it was shown that the complex of curves, $\mathcal{C}(S)$, is $\delta$-hyperbolic. In contrast, in section \ref{Examples} it will be shown that

\begin{thm}For $g>3$ and $\alpha\neq 0$, $\mathcal{HC}(S,\alpha)$ is not $\delta$-hyperbolic. 
\label{thmone}
\end{thm}

It is also well known (for example \cite{MasurandMinskyI}) that in $\mathcal{C}(S)$, the distance between two vertices representing the curves $a$ and $b$ is bounded from above by the logarithm of the intersection number (see section ``Intersection numbers'' for definition). However, in section \ref{Examples} it will be shown that

\begin{thm}Let $m_1$ and $m_2$ be multicurves in the integral homology class $\alpha$. Then $d_{\mathcal{H}}(m_{1},m_{2})\leq \frac{i(m_{1},m_{2})}{2}+1$, where $i(m_{1},m_{2})$ is the geometric intersection number. This bound is sharp.
\label{thmtwo}
\end{thm}

\begin{figure}[h]
\centering
\includegraphics[width=13cm]{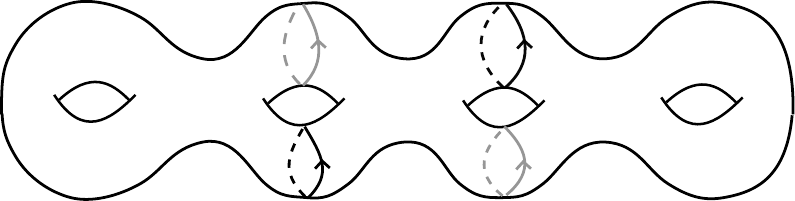}
\caption{An edge that is not simple. The multicurve drawn in grey represents one vertex and the multicurve drawn in black represents the other.}
\label{chone1}
\end{figure}

An edge in $\mathcal{HC}(S,\alpha)$ connecting two vertices representing the multicurves $\gamma_{i}$ and $\gamma_{i+1}$ is called \textit{simple} if $\gamma_{i+1}-\gamma_i$ is the oriented boundary of an embedded subsurface of $S$, see figure \ref{chone1}. A \textit{simple path} is a path that only traverses simple edges. In section \ref{algorithm} an algorithm for constructing simple paths between any two vertices (hereafter referred to as the ``path construction algorithm'') is given. 

Let $I$ be a closed interval. In section \ref{surfaceconstruction} a path in $\mathcal{HC}(S,\alpha)$ connecting the vertices representing $m_1$ and $m_2$ is shown to correspond to an oriented, immersed surface $H$ in $S\times I$ with $\partial H$ homotopic to the multicurves $m_{2}-m_{1}$ in $S\times 0$. The geometry of $\mathcal{HC}(S,\alpha)$ is thus related to the topology of surfaces in $S\times I$. 

In \cite{Me2}, it is shown that every oriented, embedded, incompressible surface in $S\times \mathbb{R}$ with boundary $m_{2}-m_{1}$ can be constructed from a path in $\mathcal{HC}(S,\alpha)$. In section \ref{surfacesection}, the path construction algorithm is used to prove that

\begin{thm}
Consider the set of all homotopy classes of orientable surfaces in $S\times \mathbb{R}$ with boundary $m_{2}-m_{1}$. Let $\mathcal{F}$ be the subset with minimal genus. Then $\mathcal{F}$ always contains an embedded surface.
\label{afterthought}
\end{thm}

Modulo a uniformly bounded multiplicative constant, it follows that the distance between two vertices in $\mathcal{HC}(S, \alpha)$ representing the multicurves $m_1$ and $m_2$ is equal to the smallest possible genus of an orientable surface in $S\times I$ with boundary $m_{2}-m_{1}$. An explicit algorithm for constructing the embedded, minimal genus surface from theorem \ref{afterthought} is given in section \ref{surfacesection}.\\

In order to show that the path construction algorithm is optimal in some sense, the geometry of $\mathcal{HC}(S, \alpha)$ is related to the topology of immersed surfaces in $S\times I$ by defining two functions from $S\setminus (m_{1}\cup m_{2})\rightarrow \mathbb{Z}$: the overlap function and the pre-image function. These functions will now be briefly described.\\

\textbf{Intersection numbers}. There are two types of intersection numbers used in this work. The \textit{intersection number}, also known as the geometric intersection number, of two multicurves $m_1$ and $m_2$ is the minimum possible number of intersections between a pair of multicurves, one of which is isotopic to $m_1$ and the other to $m_2$. The intersection number of $m_1$ and $m_2$ is denoted by $i(m_{1},m_{2})$, and the \textit{algebraic intersection number} is denoted by $\hat{i}(m_{1},m_{2})$. The algebraic intersection number of an oriented arc $a$ with an oriented representative $m_1$ of the homotopy class $[m_{1}]$ is also written as $\hat{i}(a,m_{1})$. \\

Intersection numbers of curves in $S\times I$ are defined by projecting onto $S\times 0$. A union of cycles in $S\times I$ is defined to be a multicurve if it projects onto a multicurve in $S\times 0$.\\

\textbf{The pre-image function.} Let $\pi$ be the projection of $S\times I$ onto $S\times 0$ given by $(s,r) \mapsto s\times 0$. Informally, given an oriented, immersed surface $H$ in $S\times I$, the \textit{pre-image function}, $g_{H}: S\times 0 \setminus \pi(\partial H)\rightarrow \mathbb{Z}$ is given by $g_{H}(s)=\hat{i}(\pi^{-1}(s), H)$ (see Section \ref{mainargument} for a more precise definition). It is shown that, modulo an additive constant, the pre-image function does not depend on $H$ but only on its boundary (lemma \ref{function}).\\

\textbf{The overlap function, and the homological distance.} The \textit{overlap function}, also denoted by the symbol $f$, of a null homologous union of curves, $n$, is a locally constant function defined on $S\setminus n$ with minimum value zero. For any two points $x$ and $y$ in $S\setminus n$, $f(x)-f(y)$ is the algebraic intersection number of $n$ with an oriented arc with starting point $y$ and endpoint $x$. An important special case is the overlap function of the difference of two homologous multicurves, $m_{1}$ and $m_2$. \\

The overlap function is not dependent on the choice of oriented arc, because the algebraic intersection number of any closed loop with $n$ is zero. It does however depend on the choice of representatives of the homotopy classes of curves. It will be assumed that the representatives of the homotopy classes are chosen so that the maximum, $M$, of the overlap function is as small as possible. For two homologous multicurves $m_1$ and $m_2$, the quantity $M$ will be called the \textit{homological distance}, $\delta(m_{1}, m_{2})$, between $m_1$ and $m_2$. \\

If $H$ is a surface constructed from a simple path connecting $m_1$ and $m_2$, as described in subsection \ref{surfaceconstruction}, the relation between the pre-image function and the overlap function of $m_1$ and $m_2$ is used to show that the path construction algorithm constructs the shortest possible simple paths. \\

\begin{thm}
Let $m_{1}$ and $m_{2}$ be two multicurves corresponding to vertices of $\mathcal{HC}(S, \alpha)$. The shortest simple path connecting the vertices has length equal to $\delta(m_{1}, m_{2})$. 
\label{j}
\end{thm}

The path construction algorithm is similar to a construction in \cite{Hatcher} for showing contractibility of the cyclic cycle complex, and can also be used to construct paths in this complex. It will be shown in the appendix that the paths so constructed in the cyclic cycle complex are geodesics.\\

A nice property of the path construction algorithm is that, as shown in theorem \ref{backtofront}, it constructs the same unoriented path from $m_1$ to $m_2$ as from $m_2$ to $m_1$.

One reason for being interested in simple paths is that they give a good estimate of distance.\\

\begin{thm}

Suppose $m_1$ and $m_2$ are homologous multicurves neither of which contain null homologous submulticurves or homotopic curves. Then $\frac{1}{-3\chi(S)}\delta(m_{1}, m_{2})\leq d(m_{1}, m_{2})\leq \delta(m_{1}, m_{2})$.
\label{quasigeodesic}
\end{thm}



\textbf{The Case} $\alpha=0$. The case in which $\alpha$ is allowed to be zero is quite different. For example, in this case the complex admits an action of the full mapping class group, and when $\alpha$ is nontrivial, it does not. In the latter case, the natural group that acts is the subgroup of the mapping class group preserving $\alpha$. 
Various complexes of null homologous (multi)curves, have been studied, for example the complex of separating curves and the Torelli geometry. Some of the methods discussed in this paper generalise, however the main problem seems to be that performing surgeries on null homologous multicurves could give null homologous curves.

\subsection*{Acknowledgements} I would like to thank Ursula Hamenst\"adt for her supervision of this project. Also, without the advice and enthusiasm of many people, the writing of this paper could have dragged on into infinity. Thanks to Joan Birman, Carl-Friedrich B\"odigheimer, Benson Farb, Sebastian Hensel, Andrew Putman and Kasra Rafi. I am particularly grateful to Dan Margalit for his patience in teaching me how to write papers, and to Allen Hatcher and an anonymous reviewer for their detailed comments and improvements.

\section{Simple Paths}
In this section, the notion of a simple path is introduced in order to be able to perform counting arguments that relate surfaces in $S\times I$ to paths in $\mathcal{HC}(S,\alpha)$.\\

A \textit{curve} $c$ in $S$ is a piecewise smooth, injective map of $S^1$ into $S$ that is not null homotopic. A \textit{multicurve} is a union of pairwise disjoint curves on $S$, and is allowed to contain null homologous submulticurves. When convenient, a curve is confused with its image in $S$.\\

Whenever this does not lead to confusion, the same symbol will be used for a vertex in $\mathcal{HC}(S,\alpha)$ and the corresponding multicurve on $S$. Also, a path in $\mathcal{HC}(S,\alpha)$ will often be denoted by a sequence of multicurves, $m_{1}, m_{2},\ldots, m_{n}$ with the property that $m_i$ and $m_{i+1}$ are disjoint for every $1\leq i \leq n$, i.e. $m_i$ and $m_{i+1}$ represent an edge in $\mathcal{HC}(S, \alpha)$.\\

If a null homologous multicurve $n$ bounds an embedded subsurface of $S$, the union of the components of $S\setminus n$ whose boundary orientation coincides with the orientation of $n$ will be called the \textit{subsurface of $S$ bounded by $n$.} If $n$ contains homotopic curves with opposite orientation, these curves are thought of as bounding an annulus, not the empty set. This convention ensures that surfaces constructed from simple paths, as outlined in section \ref{surfaceconstruction}, are embedded.\\

The next lemma is used to decompose null homologous multicurves into boundaries of subsurfaces. 

\begin{lem}
\label{boundary}
If a null homologous multicurve $n$ does not contain a nontrivial null homologous submulticurve, it bounds a subsurface of $S$.
\end{lem}

\begin{proof}
Consider the subsurface of $S$ on which the overlap function of $n$ has its maximum. Its boundary is a null homologous submulticurve of $n$. By assumption on $n$ it must be all of $n$.


\end{proof}

\begin{cor}
For any path $v_{1}, v_{2},\ldots, v_n$ in $\mathcal{HC}(S,\alpha)$, a simple path can be obtained by adding extra vertices where necessary.
\end{cor}

\begin{proof}
Suppose $m_1$ and $m_2$ are connected by an edge that is not simple. By the previous lemma, $m_{2}-m_{1}$ can be decomposed into $k$ null homologous submulticurves $n_{1}, n_{2},\ldots,n_{k}$, each of which bounds a subsurface of $S$. Then a simple path connecting $m_1$ and $m_2$ is determined by the vertices $m_{1}, m_{1}+n_{1}, m_{1}+n_{1} +n_{2},\ldots,m_{1}+n_{1}+n_{2}+\ldots+n_{k-1},m_{2}$.
\end{proof}



\subsection{Constructing an Embedded Surface in $S\times I$ from a Path in $\mathcal{HC}(S,\alpha)$}
\label{surfaceconstruction}
All curves, surfaces, and manifolds discussed here are assumed to be piecewise smooth.



Suppose $\gamma$ is a simple path in $\mathcal{HC}(S,\alpha)$ passing through the vertices corresponding to the multicurves $\gamma_{0}, \gamma_{1}, \ldots,\gamma_{j}$. A surface $T_{\gamma}$ contained in $S\times j$ is constructed inductively. Given $\gamma_0$, isotope $\gamma_1$ such that there is a subsurface $S_1$ of $S$ with boundary $\gamma_{1}-\gamma_0$. Let $T_1$ be the surface in $S\times [0,1]$ given by $\gamma_{0}\times [0,\frac{1}{2}]\cup S_{1}\times \{\frac{1}{2}\}\cup \gamma_{1}\times [\frac{1}{2},1]$. Next, isotope $\gamma_2$ so that there is a subsurface $S_2$ of $S$ with $\partial S_{2}=\gamma_{2}-\gamma_1$ and let $T_{2}=\gamma_{1}\times [1,\frac{3}{2}]\cup S_{2}\times \{\frac{3}{2}\}\cup \gamma_{2}\times [\frac{3}{2},2]$. Repeat this successively for each of the $\gamma_i$ until an embedded surface $T_{\gamma}=T_{1}\cup T_{2}\cup\ldots\cup T_{j}$ in $S\times [0,j]$ is obtained.\\

$T_{\gamma}$ is called the \textit{trace surface} of the path $\gamma$. Note that the trace surface of a path depends on the orientation on $S$.\\


\textbf{Remark} Similarly, if $\gamma_{0}, \gamma_{1},\ldots,\gamma_j$ is not simple, the above procedure can be used to construct a cell complex with boundary $\gamma_{j}-\gamma_{0}$. It is not difficult to show that such cell complexes are homotopic to immersed surfaces in $S\times [0,j]$.

\subsection{Extrema of the overlap function}
In order to construct paths in $\mathcal{HC}(S, \alpha)$, it is necessary to use some properties of the level sets, in particular the local extrema, of the overlap function. These will be used to define the surgeries used in the path construction algorithm.\\

Given an oriented multicurve $a$ with a regular neighbourhood $\mathcal{N}(a)$ and an orientation on $S$, the left and right component of $\mathcal{N}\setminus a$ can be defined. If $b$ is an oriented multicurve that intersects $a$ transversely at a point $p$, it therefore makes sense to say that $b$ \textit{crosses over} $a$ \textit{from left to right} (or \textit{right to left}) at $p$. Similarly, if $b$ is an oriented arc with an endpoint on $a$, a notion in which $b$ \textit{leaves or approaches} $a$ \textit{from the left or right} can be defined.\\

If a horizontal arc of $a\cap (S\setminus b)$ leaves and approaches $b$ from the right, then this arc is \textit{to the right of }$b$ and vice versa. \\

Whenever $m_1$ and $m_2$ are homologous multicurves, the overlap function of $m_1$ and $m_2$ is bounded and has a maximum. Call the subsurface of $S$ on which the overlap function takes on its maximum $S_{max}(m_{1},m_{2})$. $S_{max}(m_{1},m_{2})$ has at least one connected component. The boundary of $S_{max}(m_{1},m_{2})$ consists of arcs of $m_1$ and $m_2$ such that $S_{max}(m_{1},m_{2})$ is to the right of any arc of $m_1$ on its boundary and to the left of any arc of $m_2$ on its boundary. In other words, the boundary of $S_{max}(m_{1},m_{2})$ is a null homologous multicurve made up of arcs of $m_1$ to the left of $m_2$ and arcs of $m_2$ to the right of $m_1$. \\

\begin{figure}
\begin{center}
\def\svgwidth{13cm}
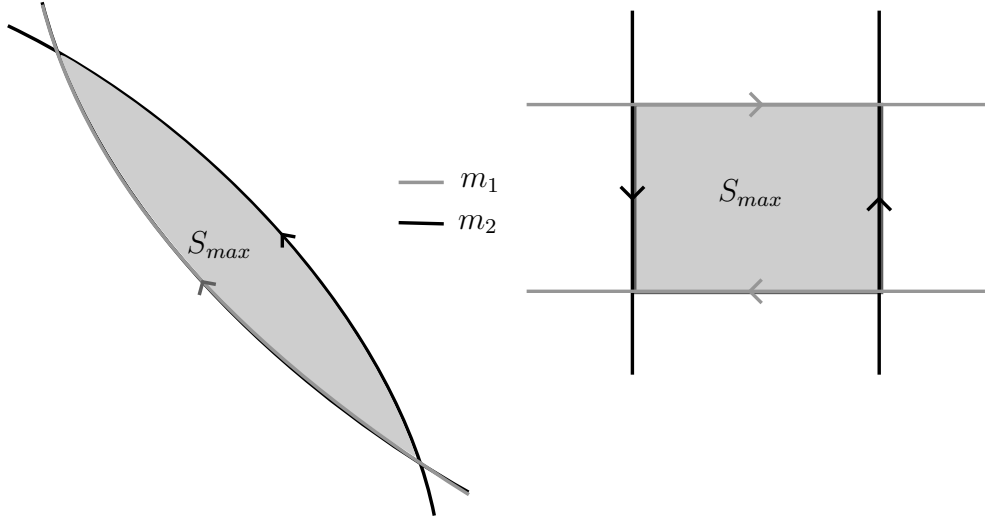
\caption{Examples of $S_{max}(m_{1},m_{2})$.}
\label{chone6}
\end{center}
\end{figure}

Similarly, the subsurface of $S$, $S_{min}(m_{1},m_{2})$, on which $f=0$, is disjoint from $S_{max}(m_{1},m_{2})$ and is on the left of any arc of $m_1$ on its boundary and to the right of any arc of $m_2$ on its boundary.

\subsection{Horizontal and vertical arcs}
Given two multicurves $a$ and $b$ on an oriented surface $S$, a \textit{horizontal arc} of $a$ is a component of $a\cap (S\setminus b)$ that leaves and approaches $b$ from the same side. A \textit{vertical arc} of $a\cap (S\setminus b)$ leaves and approaches $b$ from opposite sides. An ``innermost'' arc in \cite{Hatcher} is an example of a horizontal arc.\\

\begin{figure}
\begin{center}
\def\svgwidth{13cm}
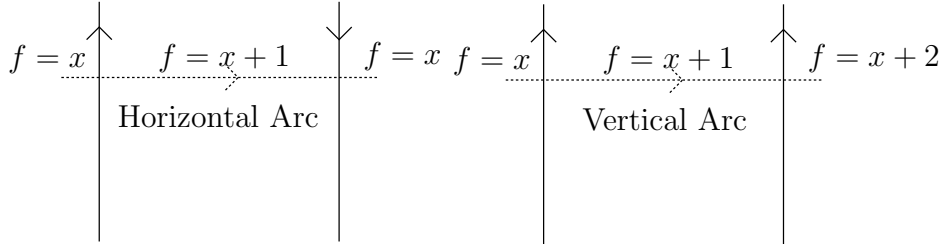
\caption{When the overlap function is thought of as a height function, a horizontal arc is horizontal with respect to this height function, and a vertical arc is vertical.}
\label{horizontalvertical}
\end{center}
\end{figure}

Suppose $a$ and $b$ are multicurves in $S$ in general position. Two arcs $a_1$ and $a_2$ of $a\cap (S\setminus b)$ will be called homotopic if the closure of $a_1$, $\bar{a}_1$, can be homotoped onto the closure of $a_2$, $\bar{a}_2$, by a homotopy that keeps the endpoints of the arcs on $b$. Since an arc of $a$ is defined to be a connected component of $a\cap (S\setminus b)$, a homotopy is also not allowed to move any interior point of the arc over $b$. Two oriented arcs will be said to be \textit{homotopic and oriented in the same way} if one can be homotoped into the other in such a way that the orientations coincide.\\

It is not difficult to see that the property of being horizontal or vertical is invariant under homotopy. Also, a horizontal arc of $a\cap (S\setminus b)$ to the right of $b$ can not be homotopic to a horizontal arc of $a\cap (S\setminus b)$ to the left of $b$, and an oriented arc of $a\cap (S\setminus b)$ is not homotopic to itself with the opposite orientation.\\

The arcs on $\partial S_{max}(m_{1},m_{2})$ and $\partial S_{min}(m_{1},m_{2})$ are all horizontal.

\subsection{Minimising Overlap}
A difficulty is that vertices of $\mathcal{HC}(S,\alpha)$ are only defined up to homotopy, whereas some of the quantities, such as the overlap function, also depend on the representative of the homotopy classes. For this reason it is necessary to work with representatives of the free homotopy class that minimise the overlap function. 

\begin{figure}
\begin{center}
\def\svgwidth{13cm}
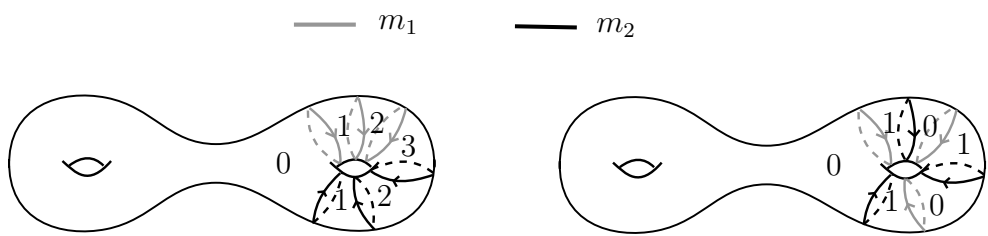
\caption{On the left, $m_1$ and $m_2$ are not in minimal position, because the overlap function could be made smaller, as shown on the right. The numbers shown are the values of the overlap function.}
\label{chone5}
\end{center}
\end{figure}

Two multicurves $m_1$ and $m_2$ will be said to be in \textit{minimal position} if
\begin{itemize}
\item $m_1$ and $m_2$ are in general position
\item the number of times $m_1$ intersects $m_2$ is equal to $i(m_{1},m_{2})$, and
\item whenever $m_1$ and $m_2$ are homolgous and $m_{2}-m_{1}$ contains homotopic curves, these homotopic curves are positioned in such a way that the overlap function is minimised. An example is illustrated in figure \ref{chone5}.
\end{itemize}

\section{A path constructing algorithm}
\label{algorithm}

\begin{figure}[h]
\begin{center}
\def\svgwidth{13cm}
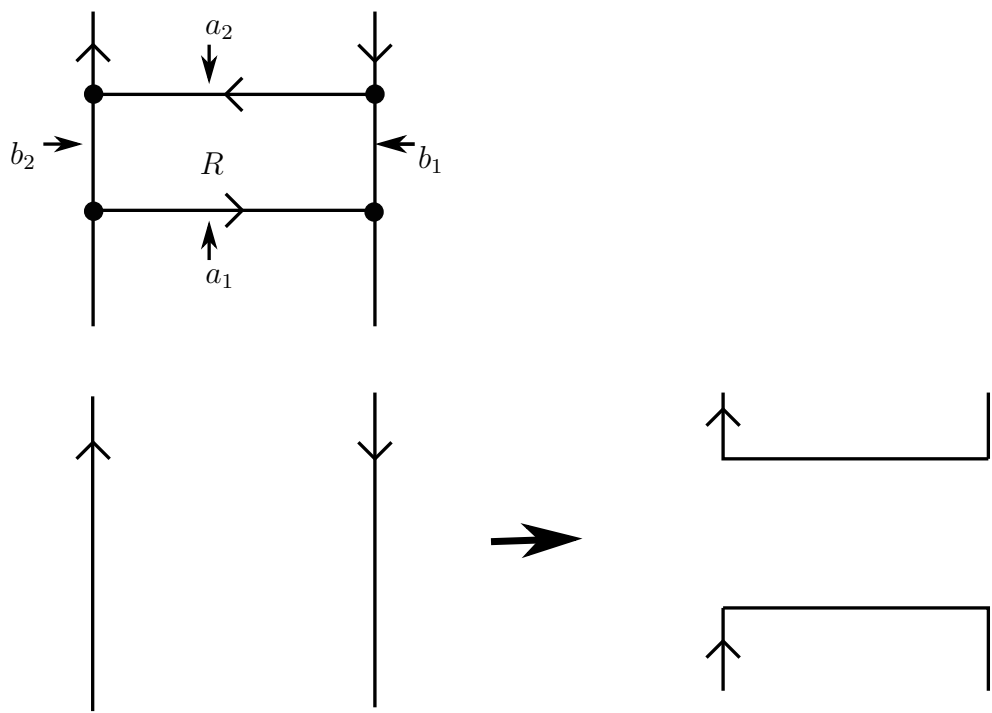
\caption{Surgering a multicurve along a horizontal arc.}
\label{chone4}
\end{center}
\end{figure}

In this section an algorithm for constructing a simple path $m_{1}$, $\gamma_{1},$ $\gamma_{2}$,\ldots,$m_{2}$ of length $\delta(m_{1}, m_{2})$ will be constructed. Recall from the introduction that $\delta(m_{1}, m_{2})$ is equal to the maximum of the overlap function of $m_1$ and $m_2$.\\

\textbf{A basic surgery construction.} Suppose $b$ is an oriented multicurve, and $a_1$ and $a_2$ are two homotopic arcs with endpoints on $b$. Then $a_{1}\cup a_{2}\cup b$ can be thought of as a one dimensional simplicial complex on $S$. If $a_1$ is a horizontal arc, the arcs $a_1$ and $a_2$ can be oriented in such a way that $a_{1}\cup a_{2}\cup b_{1}\cup b_{2}$ bound a rectangle $R$ in $S$, where $b_1$ and $b_2$ are chains in $a_{1}\cup a_{2}\cup b$, as shown in figure \ref{chone4}. \textit{Surgering an oriented multicurve} $b$ \textit{along a horizontal arc} $a$ is the process in which the oriented chains $-a_1$, $-a_2$, $-b_1$ and $-b_2$ are added to the subcomplex $b$. Since the chain added is a boundary, the resulting multicurve is homologous to $b$. 
\\


Suppose $m_1$ and $m_2$ are multicurves in minimal position. The union of multicurves, $m_{2}-m_{1}$, defines a one dimensional cell complex on $S$. By convention, $\partial S_{max}(m_{1},m_{2})$ is oriented in such a way that $S_{max}(m_{1},m_{2})$ is on its left. Let $a_1, a_2,\ldots$ be the arcs of $m_2$ on $\partial S_{max}(m_{1},m_{2})$, and $b_1, b_2, \ldots$ be the arcs of $m_1$ on $\partial S_{max}(m_{1},m_{2})$. Then $\partial S_{max}(m_{1},m_{2})= \sum_i a_i - \sum_j b_j$ is a union of chains. The multicurve $\gamma_1$ is obtained by adding $\partial S_{max}(m_{1},m_{2})$ to $m_1$ as a chain, i.e. $\gamma_{1}$ is the subcomplex $m_{1}+\cup a_{i}-\cup b_{i}$. The surgery in which $\gamma_1$ is constructed from $m_1$ and $m_2$ is called \textit{performing the surgery} or \textit{surgeries} corresponding to $S_{max}(m_{1},m_{2})$ on $m_1$. 

An equivalent means of constructing $\gamma_1$ is as follows:
\begin{enumerate}
\item The multicurve $m_1$ is first surgered along the arcs $a_1, a_{2},\ldots$. 
\item The surgery from one gives a multicurve $\gamma_{1}-\partial S_{max}(m_{1},m_{2})$. Discard the null homologous submulticurve $-\partial S_{max}(m_{1},m_{2})$.
\end{enumerate}

Up to free homotopy on the boundary, $S_{max}(m_{1},m_{2})$ can be thought of as ``that piece of $S$ that is bounded by $m_1$ and $\gamma_1$''.

\begin{figure}
\begin{center}
\def\svgwidth{10cm}
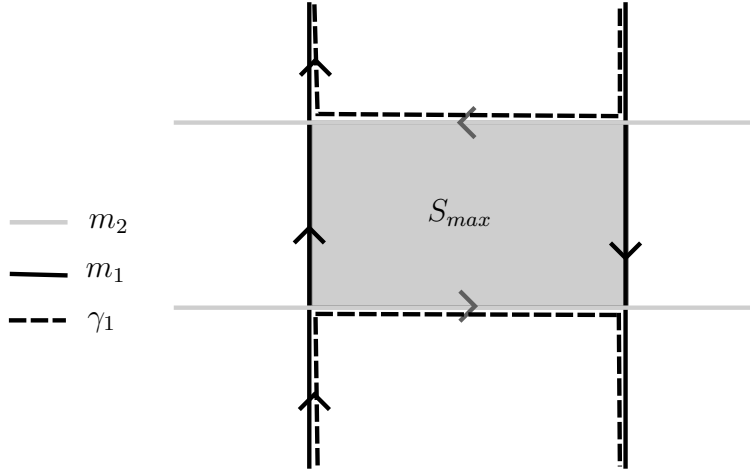
\caption{How to construct $\gamma_{1}$}
\label{chone7}
\end{center}
\end{figure}


By construction, $i(\partial S_{max}(m_{1},m_{2}),m_{1})=0$ and each connected component of $S_{max}(m_{1},m_{2})$ intersects an annular neighbourhood of $m_1$ on the right side of $m_1$ (i.e. every component of $S_{max}(m_{1},m_{2})$ is ``on the same side'' of $m_1$). Therefore $i(\gamma_1,m_1)=0$. \\

As constructed, the multicurve $\gamma_{1}$ might contain trivial curves that bound disks, and might not be in minimal position with $m_2$. This point is ignored at the moment. Only once all the multicurves $\gamma_i$ are constructed are the trivial curves discarded from each $\gamma_i$. \\

Cutting out the arcs $b_i$ make it possible to connect the subsurface of $S$, $S_{min}(\gamma_{1},m_{2})$, on which $f_1$ takes on its minimum, to $S_{max}(\gamma_{1},m_{2})$ (defined similarly), by an arc that crosses $m_2-\gamma_1 $ from right to left once less than any arc connecting $S_{min}(m_{1},m_{2})$ with $S_{max}(m_{1},m_{2})$. In other words, $\delta(\gamma_{1}, m_{2})=\delta(m_{1},m_{2})-1$.\\

Let $f_k$ be the overlap function of $\gamma_k$ and $m_2$. The multicurve $\gamma_{i+1}$ is constructed from $\gamma_i$ in the same way as $\gamma_1$ from $m_1$ only with the multicurve $m_1$ replaced by $\gamma_i$. The \textit{surgery corresponding to} $S_{max}(\gamma_{i},m_{2})$ is defined to be the surgery performed on $\gamma_{i}$ to obtain $\gamma_{i+1}$.\\

\begin{figure}[h]
\begin{center}
\def\svgwidth{11cm}
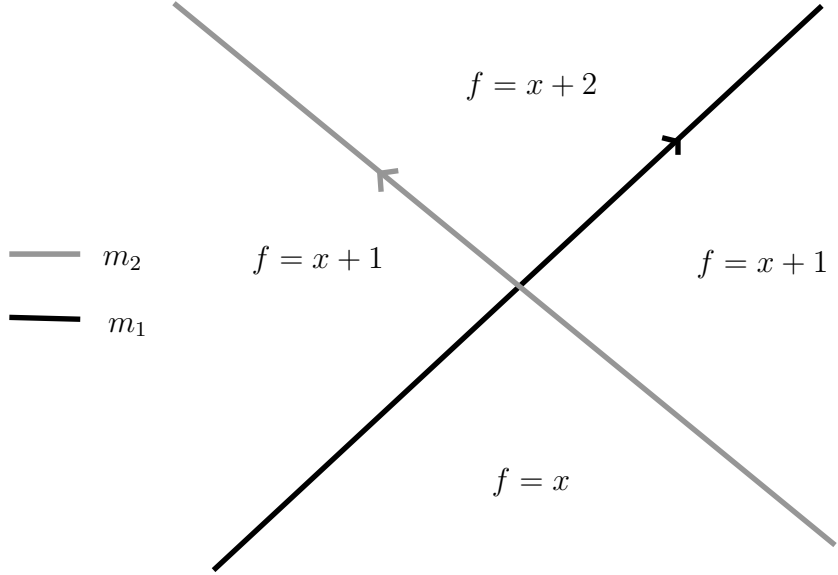
\caption{A point of intersection forces the overlap to have maximum at least two.}
\label{chone8}
\end{center}
\end{figure}

The construction ends with the multicurve $\gamma_j$ when $\delta(\gamma_{j},m_{2})=1$. This can only be happen if $\gamma_j$ and $m_2$ do not intersect, because as shown in figure \ref{chone8}, an intersection forces the maximum of $f_j$ to be at least two.\\

If $\delta(\gamma_{j},m_{2})=1$, then $S_{max}(\gamma_{j},m_{2})$ is the subsurface bounded by $m_2-\gamma_j$. \\

\textbf{Remark} The reason for giving two equivalent definitions of surgery corresponding to $S_{max}$ is to make it clear that for any two homologous multicurves $m_1$ and $m_2$, a path between $m_1$ and $m_2$ in $\mathcal{HC}(S,\alpha)$ can be constructed by repeatedly surgering along horizontal arcs and adding/discarding null homologous submulticurves. Since this also applies to multicurves representing vertices joined by an edge in $\mathcal{HC}(S,\alpha)$, it follows that every path in $\mathcal{HC}(S,\alpha)$ can be constructed by surgering along arcs and adding or discarding null homologous submulticurves. This is the approach taken in \cite{Hatcher}, and will be used in the proof of theorem \ref{quasigeodesic}.\\

\textbf{Remark} The choice to use $S_{max}(m_{1},m_{2})$ instead of $S_{min}(m_{1},m_{2})$ was arbitrary. However, it is not possible to simultaneously reduce the intersection number further at each step by requiring that the subsurface of $S$ bounded by $\gamma_1$ and $m_1$ be $S_{max}(m_{1},m_{2})\cup S_{min}(m_{1},m_{2})$. This is because $S_{min}(m_{1},m_{2})$ is to the left of $m_1$ and $S_{max}(m_{1},m_{2})$ is to the right of $m_1$, so a simple path would not be obtained.\\

This completes the construction of the promised algorithm. A simple path constructed in this way will be called a \textit{middle path}. The algorithm itself will be referred to as the \textit{path construction algorithm}.\\

\begin{figure}
\begin{center}
\def\svgwidth{12cm}
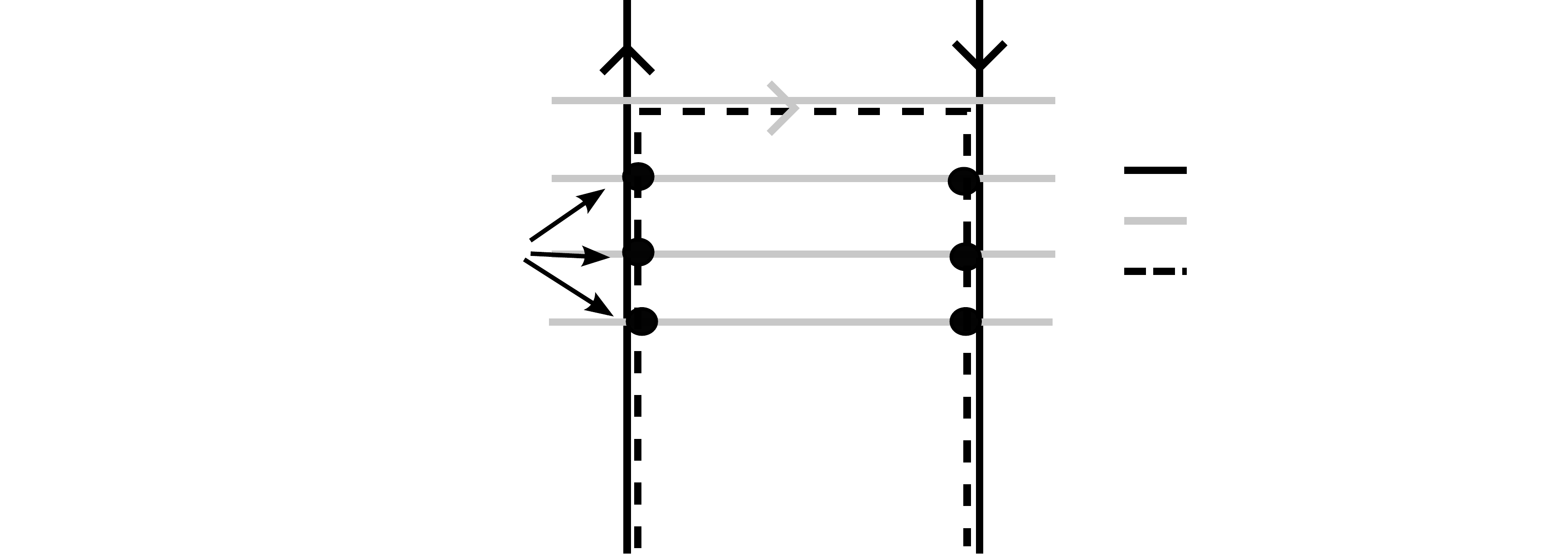
\label{chone11redone}
\end{center}
\end{figure}

\textbf{Minimal position problems.} The multicurve $\gamma_{i+1}$, obtained from $\gamma_i$ and $m_{2}$ by performing the surgery corresponding to $S_{max}(\gamma_{i},m_{2})$, might contain curves that bound disks or there might be points of intersection with $m_{2}$ that can be removed by a homotopy, as shown in figure \ref{chone11redone}. In other words, the multicurves $\gamma_i$ are not in minimal position. 

Sometimes it is convenient to drop the assumption that multicurves are in minimal position, and only require that the multicurves are chosen so as to minimise the maximum of the overlap function. For example, in the proof of theorem \ref{backtofront}, it is desirable to have $\gamma_{i}$ a subcomplex of the one dimensional cell complex $\gamma_{1}\cup \gamma_{k}$. This choice of $\left\{\gamma_{i}\right\}$ are not in general position and can have points of intersection that can be removed by a homotopy.\\

We assume that $m_1$ and $m_2$ are in minimal position, and show that the multicurves constructed by the path construction algorithm are representatives of their homotopy class that minimise the maximum of the overlap function with $m_1$ and $m_2$. Performing the surgery corresponding to $S_{max}(m_{1},m_{2})$ on $m_1$ gives the multicurve $\gamma_{1}$. The overlap function calculated from $\gamma_1$ and $m_2$ has maximum equal to $\delta(m_{1},m_{2})-1$. There is no multicurve $\gamma_{1}^{'}$ in the homotopy class $[\gamma_{1}]$ with the property that the overlap function of $\gamma_{1}^{'}$ with $m_2$ has maximum less than $\delta(m_{1},m_{2})-1$. For the maximum of $\gamma_{1}^{'}$ to be less than $\delta(m_{1},m_{2})-1$, $\gamma_{1}^{'}-m_{1}$ would have to be a null homologous multicurve $n$ with the property that every arc connecting $S_{max}(m_{1},m_{2})$ with $S_{min}(m_{1},m_{2})$ intersects $n$ at least twice. This is not possible because $\gamma_{1}^{'}$ is homotopic to $\gamma_{1}$, and $\gamma_{1}-m_{1}$ is the boundary of an embedded subsurface of $S$. It follows from the same argument that the maximum of the overlap function of $\gamma_2$ and $m_2$ is equal to $\delta(\gamma_{2}, m_{2})$, similarly for $\gamma_{3}$ and $m_2$, etc.\\

The next lemma is needed in theorems in which it is necessary to compare the overlap functions $f_i$ for different values of $i$.

\begin{lem}
Let $\gamma_{1}, \gamma_{2},\ldots,\gamma_{k}$ be a middle path. There is a representative $\gamma_i$ of the homotopy class $[\gamma_{i}]$ with the property that $\gamma_i$ is an oriented, embedded subcomplex of the one dimensional oriented cell complex $\gamma_{1}\cup \gamma_{k}$.
\label{gp}
\end{lem}

\begin{proof}
Suppose $\gamma_1$ and $\gamma_k$ are in minimal position. Choose the representatives $\gamma_{i+1}$ and $\gamma_i$ of $[\gamma_{i+1}]$ and $[\gamma_{i}]$ such that $\gamma_{i+1}-\gamma_{i}= \partial S_{max}(\gamma_{i},m_{2})$. The boundary of $S_{max}(\gamma_{i},m_{2})$ is an embedded subcomplex of $\gamma_{1} \cup \gamma_k$ for every $i$, and has zero intersection number with $\gamma_{i}$ and $\gamma_{i+1}$. Recall that the multicurve $\gamma_{i+1}$ is obtained from $\gamma_i$ by subtracting the arcs of $\gamma_{i}\cap (S\setminus \gamma_{k})$ on $\partial S_{max}(\gamma_{i},m_{2})$ and adding the arcs of $\gamma_{k}\cap (S\setminus \gamma_{i})$ on $\partial S_{max}(\gamma_{i},m_{2})$. Also, no arc of $\gamma_{k}\cap \gamma_{1}$ will be on the boundary of $S_{max}(\gamma_{i},m_{2})$ for more than one $i$, so each arc can only be added or subtracted at most once. Each of the multicurves $\gamma_{i}$ is therefore an oriented subcomplex of $\gamma_{1}\cup \gamma_{k}$. From figure \ref{chone8}, it is easy to verify that $S_{max}(\gamma_{i},m_{2})$ can not meet itself at a vertex, because if four components of $S\setminus (\gamma_{i}\cup \gamma_{k})$ come together at a point and the overlap function is equal on two of them, it must be larger on a third component and smaller on the fourth. Therefore, if $\gamma_i$ doesn't meet or cross over itself at a vertex, neither will $\gamma_{i+1}$. The $\gamma_i$ chosen in this way are therefore also embedded. 

\end{proof}



A nice property of the path construction algorithm is that it constructs the same path in reverse.
\begin{thm}
\label{backtofront}
If $m_1$ and $m_2$ had been interchanged in the path construction algorithm, the same unoriented path would have been obtained.
\end{thm}

\begin{proof}
Suppose the representatives $m_{1}, \gamma_{1},\ldots,\gamma_{j}, m_{2}$ of the free homotopy classes $[m_{1}], [\gamma_{1}],\ldots,[\gamma_{j}], [m_{2}]$ are chosen as outlined in lemma \ref{gp}. In particular, each of the $\gamma_i$ are oriented subcomplexes of the cell complex $m_{1}\cup m_{2}$ such that $\gamma_{i+1}-\gamma_{i}$ is the boundary of the subsurface of $S$ on which the overlap function $f$ of $m_{2}-m_{1}$ is no less than its maximum value minus $i$.
Let $h$ be the overlap function of $m_{1}-m_{2}$. It is easy to check that $h$ has its maximum where the overlap function of $m_{2}-m_{1}$ has its minimum, and vice versa. By definition, $\gamma_j$ is the multicurve chosen such that $m_{2}-\gamma_j$ bounds the subsurface of $S$ given by $S\setminus S_{min}(m_{1},m_{2})$. In other words, $\gamma_j - m_{2}$ is the boundary of $S_{min}(m_{1},m_{2})$ or the boundary of the subsurface of $S$ on which $h$ has its maximum. The multicurve $\gamma_j$ therefore satisfies the definition of the first multicurve in the path $m_{2},\ldots,m_{1}$. Similarly for $\gamma_{j-1}$, $\gamma_{j-2}$, etc.

\end{proof}

\section{The Overlap Function and the Pre-image Function}
\label{mainargument}
Let $H$ be an oriented, immersed surface in $S\times I$, where $I$ is the interval $\left[0, N\right]$. Suppose also that $\partial H \subset S\times \partial I$, $\pi(\partial H \cap (S\times \left\{N\right\}))=m_{2}$ and $\partial H \cap (S\times \left\{0\right\})=-m_1$, where $m_1$ and $m_2$ are homologous multicurves. \\

In this section, theorem \ref{j} is proven by relating the overlap function of $\partial H$ to the pre-image function $g_H$. \\

The \textit{pre-image function} $g_{H}:S\times 0\setminus \pi(\partial H) \rightarrow \mathbb{Z}$ is defined as follows: Suppose $P:=S\times I$ and $B$ is an open set in $(S\times \left\{0\right\})\setminus \pi(\partial H)$. Algebraic intersection number provides a map $H_{2}(P, \partial H)\times H_{1}(P,B)\to \mathbb{Z}$. For $x$ in $(S\times \left\{0\right\})\cap B$,
\begin{equation}
g_{H,B}(x):=\hat{i}(H, x\times I)
\end{equation}

For all $x\subset S\times 0 \setminus \pi(\partial H)$ there is a choice of $B$ such that $x\subset B$. The pre-image function is well defined because if $B\subset B^{'}\subset (S\times 0)\setminus \pi(\partial H)$, it follows from the naturality of the intersection pairing with respect to inclusions (\cite{Dold} Proposition 1.3.4) that the diagram below commutes.
$$
\xymatrixrowsep{1in}
\xymatrixcolsep{1in}
\xymatrix{B \ar@{^{(}->}[r] \ar[dr]^{g_{H,B}}& B^{'}\ar[d]^{g_{H,B^{'}}}\\
            & \mathbb{Z}}
$$

\begin{lem}
\label{function}
Given any two oriented, immersed surfaces $H_1$ and $H_2$ with $\pi(\partial H_{1})=\pi(\partial H_{2})=m_{2}-m_{1}$, there is a constant integer $c$ such that for all $s\in (S\times \left\{0\right\})\setminus (m_{2}-m_{1})$, we have $g_{H_1}=g_{H_2}+c$.
\end{lem}

\begin{proof}
The functions $g_{H_1}$ and $g_{H_2}$ both increase by one when crossing over an arc of $m_{2}-m_{1}$ from right to left. This lemma is thus proven by showing that $g_{H_1}$ and $g_{H_2}$ are both locally constant on $(S\times \left\{0\right\})\setminus (m_{2}-m_{1})$. Suppose $B$ is an open set in $(S\times \left\{0\right\})\setminus (m_{2}-m_{1})$ containing the points $x$ and $y$. Whenever $x$ and $y$ are points lying in the same connected component of $B$, $\{y\}\times I$ and $\{x\}\times I$ represent the same class in $H_{1}(P,B)$. It follows from the definition of $g_{H_{1}}$ that $g_{H_{1}}(x)=g_{H_{1}}(y)$, as desired. The same argument applies to $g_{H_{2}}$, from which the lemma follows.

\end{proof}

It is now possible to give a proof of theorem \ref{j}.


\begin{proof}[Proof of theorem \ref{j}]


Suppose $\gamma$ is a simple path connecting $m_1$ and $m_2$ of length less than $\delta(m_{1},m_{2})$. Let $T_{\gamma}$ be the trace surface of $\gamma$. Then $T_{\gamma}$ can be constructed by connecting up $\delta(m_{1},m_{2})-1$ or fewer pieces, each of which projects one to one onto a subsurface of $S\times 0$ with the induced subsurface orientation. Since all the subsurface glued together to form the trace surface are oriented as subsurfaces of $S\times 0$, $g_{T_{\gamma}}$ is everywhere nonnegative. It follows from lemma \ref{function} that the maximum of $g_{T_{\gamma}}$ is greater than or equal to $\delta(m_{1},m_{2})$. In other words, $\hat{i}(\pi^{-1}(s), H)\geq \delta(m_{1},m_{2})$ for some $s$. This is a contradiction.

Paths with this minimum length can always be achieved by the path construction algorithm.

\end{proof}


\section{Distances and Simple Paths}
Theorem \ref{j} determines the length of the shortest simple paths connecting two vertices, however this has not yet been related to the distance between the vertices. In order to compute distance, it is necessary to consider all (possibly nonsimple) paths.\\

A \textit{quasi-geodesic} is a map $q$ from $\mathbb{Z}\to \mathcal{HC}(S,\alpha)$ such that there are constants $K>0$ and $C>0$ such that 
\begin{equation*}
\frac{1}{K}\left|x-y\right|-C\leq d(q(x),q(y))\leq K\left|x-y\right|+C
\end{equation*}
where $d$ denotes distance in $\mathcal{HC}(S,\alpha)$. All quasi-geodesics considered in this paper are \textit{uniform quasi-geodesics}, in the sense that, for any two vertices $v_1$ and $v_2$ on the quasi-geodesic, $C=0$ and $K\leq -3\chi(S)$.\\
 
A metric space is \textit{geodesically stable}, \cite{Bonk}, if every quasi-geodesic segment is contained in the neighbourhood of a geodesic segment, where the size of the neighbourhood only depends on the constant $K$ in the definition of quasi-geodesic. Note that, since $\mathcal{HC}(S,\alpha)$ is not $\delta$-hyperbolic (in fact, it is not even nonpositively curved), no geodesic stability should be expected. Despite this, families of geodesics and/or simple paths connecting two vertices in $\mathcal{HC}(S,\alpha)$ can be easily described and constructed, however this is the subject of a future paper, \cite{Me3}.\\

Recall that if $m_{1}, \gamma_{1}, \gamma_{2},\ldots, m_{2}$ is a middle path, $\delta(\gamma_{i},m_{2})=\delta(\gamma_{i+1},m_{2})+1$. Let $m_{1}, \beta_{1}, \beta_{2},\ldots,m_{2}$ be an arbitrary (possibly nonsimple) path connecting $m_1$ and $m_2$. It is proven that middle paths are quasi-geodesics by obtaining a uniform upper bound on $\delta(\beta_{i}, m_{2})-\delta(\beta_{i+1},m_{2})$. In order to show this, the following lemma is used.\\

\begin{lem}Suppose $a$ and $b$ are multicurves in general position on $S$ such that the number of points of intersection between $a$ and $b$ is equal to $i(a,b)$. If $b$ does not contain homotopic curves, the number of homotopy classes of arcs of $a\cap (S\setminus b)$ is bounded from above by $-3\chi(S)$.
\label{beancounting}
\end{lem}
\begin{proof}
Recall that homotopy classes of arcs of $a\cap (S\setminus b)$ was defined in section 2.3. As shown in figure \ref{chfour1}, a homotopy class of arcs of $a\cap (S\setminus b)$ can be treated as a rectangle. One pair of opposite sides of the rectangle, the ``short'' sides, consist of arcs of $a\cap (S\setminus b)$ on the boundary of a component of $S\setminus (a\cup b)$ that is not a rectangle or bigon. The other pair of opposite sides, the ``long'' sides, consist of subarcs of $b$ along which the endpoints of one short side of the rectangle have to be moved by a homotopy that takes it to the opposite side of the rectangle.\\
\begin{figure}
\begin{center}
\def\svgwidth{12cm}
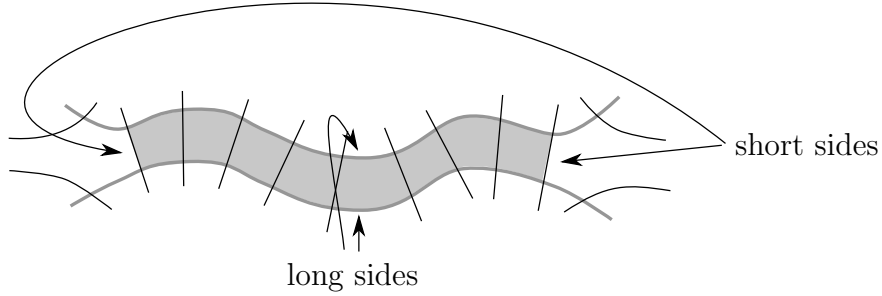
\caption{The rectangle representing a homotopy class of arcs.}
\label{chfour1}
\end{center}
\end{figure}

Each simply connected component of $S\setminus (a\cup b)$ with $2n$ sides contributes $1-\frac{n}{2}$ to the Euler characteristic of $S$; a $2n$ sided component that is not simply connected contributes even more. Apart from the rectangle, the hexagon has the largest ratio of the number of sides to its contribution to the absolute value of the Euler characteristic. The largest possible number of homotopy classes of arcs of $a\cap (S\setminus b)$ is achieved when $S\setminus (a\cup b)$ consists of rectangles and hexagons only, since every homotopy class of arcs has its short sides on the boundary of a component of $S\setminus (a\cup b)$ that is not a rectangle. In this case, the number of hexagons is equal to $-2\chi(S)$. There are three arcs of $a\cap (S\setminus b)$ on the boundary of each hexagon, and each rectangle representing a homotopy class has two short sides on the boundary of a hexagon. The bound of $-3\chi(S)$ follows directly.
\end{proof}

\begin{figure}
\begin{center}
\def\svgwidth{13cm}
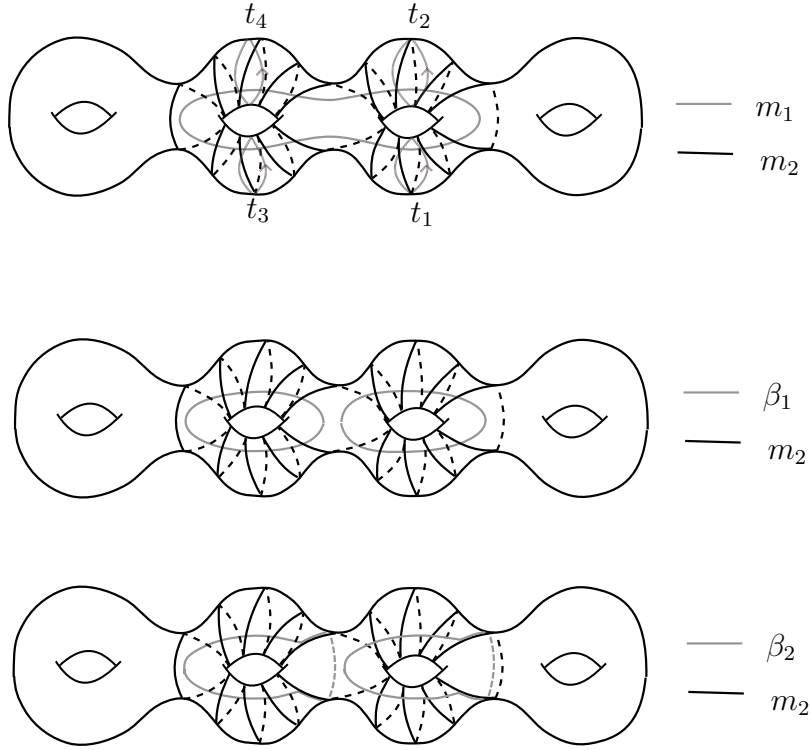
\end{center}
\label{doubletoruswithears2}
\caption{Constructing a path that is not simple}
\end{figure}

\textbf{Example - Nonsimple path} It was seen in the second last remark at the end of section 3.0 that a path in $\mathcal{HC}(S,\alpha)$ between any two vertices can be constructed by surgering along horizontal arcs and adding/discarding null homologous submulticurves. Consider the example shown in figure \ref{doubletoruswithears2}. It is possible to construct the multicurve $\beta_1$ as shown, where $\beta_1$ is not the same as the curve $\gamma_1$ constructed by the path construction algorithm. A multicurve $\beta_2$ can then be constructed by surgering $\beta_1$ along horizontal arcs, such that the number of Dehn twists of $\beta_2$ relative to $m_2$ inside each of the annuli with core curves $t_1$, $t_2$, $t_3$ and $t_4$ is one less than the number of Dehn twists of $\beta_1$ relative to $m_2$. However, the path $m_{1}, \beta_{1}, \beta_{2}$ is not simple, because $t_{1}, t_{2}, t_{3}$ and $t_4$ do not bound an embedded subsurface of $S$. In this example, the path $m_{1}, \gamma_{1},\ldots, m_{2}$ obtained from the path construction algorithm is not a geodesic. In the proof of theorem \ref{quasigeodesic}, the finite topology of the surface $S$, in the form of lemma \ref{beancounting}, is used to show that $m_{1}, \gamma_{1},\ldots, m_{2}$ is a quasi-geodesic.\\

\textbf{Performing multiple surgeries.} Suppose a multicurve $d$ is constructed from the multicurve $b$ by surgering along horizontal arcs. To be more precise, the multicurve $b$ is first surgered along a horizontal arc $a_1$ to obtain a multicurve $b_1$. The multicurve $b_1$ is then surgered along a horizontal arc $a_2$ to obtain a multicurve $b_2$, where $a_2$ is understood to be a horizontal arc in $S\setminus b_1$. The multicurve $b_2$ is then surgered along a horizontal arc $a_3$ to obtain a multicurve $b_3$, where $a_3$ is understood to be a horizontal arc in $S\setminus b_{2}$, etc, until a multicurve $b_n$ is obtained with the property that $d$ is homotopic to $b_n$. The multicurve $d$ is then referred to as \textit{the multicurve obtained by surgering $b$ along the horizontal arcs} $a_{1}, a_{2},\ldots, a_n$.

\begin{lem}
Let $\beta_{1}, \beta_{2},\ldots, \beta_n$ be a path in $\mathcal{HC}(S,\alpha)$. If $\beta_{i+1}$ can be obtained from $\beta_{i}$ by surgering along no more than $-3\chi(S)$ horizontal arcs, $a_{1}, a_{2},\ldots, a_k$, it follows that $\delta(\beta_{i},\beta_{n})-\delta(\beta_{i+1},\beta_{n})\leq-3\chi(S)$.
\label{refereelemma}
\end{lem}
\begin{proof}
Recall that the multicurve $\beta_{i+1}$ can be constructed by surgering $\beta_{i}$ along horizontal arcs, and discarding null homologous submulticurves.  Let $\beta^{'}_{i+1}$ be the multicurve obtained from $\beta_{i}$ by surgering along $a_{1}, a_{2},\ldots,a_k$. Each surgery can increase the number of curves in $\beta^{'}_{i+1}$, and hence the number of null homologous submulticurves separating $S_{max}(\beta_{i+1},m_{2})$ from $S_{min}(\beta_{i+1},m_{2})$, by no more than one. Since by assumption $\beta_{i}$ does not contain null homologous submulticurves, $\beta_{i+1}$ is obtained from $\beta_{i+1}^{'}$ by discarding no more than $-3\chi(S)$ multicurves that separate $S_{max}(\beta_{i+1},m_{2})$ from $S_{min}(\beta_{i+1},m_{2})$. The claim follows.
\end{proof}

\textbf{Remark on null homologous submulticurves and the triangle inequality.} When constructing a geodesic path $m_{1}, \beta_{1},\ldots, m_{2}$ connecting the vertices $m_1$ and $m_2$, it is possible to assume without loss of generality that the multicurves $\beta_i$ do not contain null homologous submulticurves. However, it is sometimes possible to find a union of null homologous multicurves, $N$, such that $\beta_{i}\cup N$ is a multicurve and $\delta(\beta_{i}\cup N,m_{2})<\delta(\beta_{i},m_{2})$. For this reason, when multicurves are allowed to contain null homologous submulticurves, $\delta$ does not satisfy the triangle inequality.

To construct an $N$ such that $\delta(\beta_{i}\cup N,m_{2})<\delta(\beta_{i},m_{2})$, suppose there exists a null homologous multicurve $n$ disjoint from $\beta_i$ that separates $S_{max}(\beta_{i},m_{2})$ from $S_{min}(\beta_{i},m_{2})$. $N$ consists of a union of multicurves in the homotopy class $\left[n\right]$. The orientation of the curves in $N$ is chosen such that, if $S_{max}(\beta_{i},m_{2})$ is to the right of $N$, the overlap function decreases along an arc crossing $N$ from left to right, and vice versa.

In \cite{MasurandMinskyII} the notion of a subsurface projection was defined. When $N$ is chosen to minimise $\delta(\beta_{i}\cup N,m_{2})$, $\delta(\beta_{i}\cup N,m_{2})$ is then equal to the maximum variation of the overlap function of $\beta_i$ and $m_2$ over a component of $S\setminus N$, i.e. the maximum homological distance between $\beta_i$ and $m_2$ in a subsurface projection to a component of $S\setminus N$. Since the multicurve $N$ can not contain curves from more than $-\chi(S)$ homotopy classes, it follows that $\delta(\beta_{i},m_{2})\leq -\chi(S)\delta(\beta_{i}\cup N,m_{2})$. \\


It is now finally possible to prove theorem \ref{quasigeodesic}.

\begin{proof}[Proof of theorem \ref{quasigeodesic}]
If a multicurve $m$ does not contain homotopic submulticurves, it follows from lemma \ref{beancounting} that there exists a bound of $-3\chi(S)$ on the number of pairwise disjoint homotopy classes (relative to $m$) of horizontal arcs with endpoints on $m$.\\

Suppose $m_{1}, \beta_{1},\ldots, m_{2}$ is a geodesic path connecting $m_1$ and $m_2$. Firstly, a proof of the theorem is given under the assumption that none of the $\beta_{i}$ represent multicurves containing homotopic curves. \\

Suppose $v_1$, $v_2$ are vertical arcs and $h$ is a horizontal arc, all with endpoints on $\beta_i$. Reusing the notation of lemma \ref{refereelemma}, let $a_{1}, a_{2},\ldots, a_{k}$ be arcs along which $\beta_i$ is surgered to obtain $\beta_{i+1}^{'}$. It is assumed that at least one of the arcs $a_i$ is homotopic to $v_{1}\ast h \ast v_{2}$ (otherwise the number of surgeries is automatically bounded by lemma \ref{beancounting}), and a contradiction is obtained. Lemma \ref{refereelemma} is then used to relate the number of surgeries to homological distance. \\

It is not necessary to consider trivial surgeries here, i.e. for all $j$, $S\setminus (\beta_{i}\cup a_{j})$ is not allowed to contain any bigons. For example, $\beta_{i}$ is not surgered along any two arcs in the same homotopy class. 

In $\beta_{i+1}^{'}$ there are one or two curves that were created by surgering along $v_{1}\ast h \ast v_{2}$. It can be assumed without loss of generality that $\beta_{i+1}$ contains at least one of these curves, otherwise there was no need to surger along $v_{1}\ast h \ast v_{2}$ at all.\\

Call a curve in $\beta_{i+1}$ \textit{new} if it was created by one of the surgeries in which $\beta_{i+1}$ is obtained from $\beta_i$. Either 

\begin{enumerate}
\item all new curves in $\beta_{i+1}$ are homotopic to other curves in $\beta_{i+1}$ i.e. $\beta_{i+1}$ contains homotopic curves, 
\item all new curves are homotopic to curves in $\beta_i$, i.e. $\beta_{i+1}$ is a submulticurve of $\beta_i$, or 
\item neither 1 nor 2. 

\end{enumerate}

We now show that the number of surgeries that need to be performed on $\beta_{i}$ to obtain $\beta_{i+1}^{'}$ is bounded from above by $-3\chi(S)$. Let $I$ be an oriented arc in $S$ that intersects $\beta_i$ transversely. There are a certain number of homotopy classes of arcs of $\beta_{i}\cap (S\setminus I)$ relative to $I$. Two arcs, $b_1$ and $b_2$, of $\beta_{i} \cap (S\setminus I)$ and/or $\beta_{i+1}\cap (S\setminus I)$ are defined to be homotopic if the closure of $b_1$ can be homotoped onto the closure of $b_2$ by a homotopy that keeps the endpoints on $I$. A homotopy of the multicurves $\beta_i$ or $\beta_{i+1}$ induces a homotopy of the arcs, as long as the homotopy does not take the endpoints of any arc over the boundary of $I$.\\


\begin{figure}
\begin{center}
\def\svgwidth{10cm}
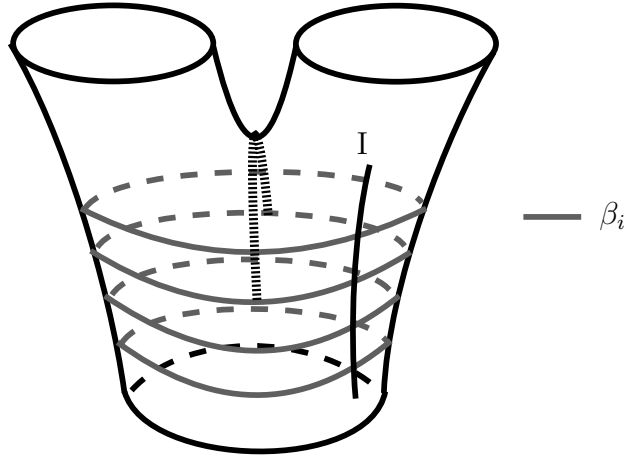
\caption{If $\beta_i$ could contain homotopic curves, the points of intersection of $\beta_{i}$ with the horizontal arc shown in the figure can be removed by a homotopy that changes the ordering of the points of intersection of $\beta_i$ with the interval $I$, without creating points of self-intersection of $\beta_i$.}
\label{chthreetwo6}
\end{center}
\end{figure}

The orientations on $I$ and $\beta_i$ make it possible to define an ordering of the starting points of the arcs of $\beta_{i}\cap (S\setminus I)$ along $I$. Suppose $I$ is chosen to contain an arc in the homotopy class $v_1$ or $v_2$ as a subarc. In the third case above, if a homotopy of $\beta_i$ or $\beta_{i+1}$ alters the order of the arcs along $I$ to remove the points of intersection with $\beta_i$ of the arc $v_{1}\ast h \ast v_{2}$, the homotopy induces points of intersection elsewhere. In other words, $i(\beta_{i+1}, \beta_{i})\neq 0$, which is not possible by definition. Since $\beta_{i+1}$ is not a submulticurve of $\beta_i$, and by assumption does not contain homotopic curves, the promised contradiction is obtained, and the claim follows in the special case that none of the $\beta_i$ contain homotopic curves. \\

As shown in figure \ref{chtwo7}, it is not always possible to get rid of all homotopic curves by assuming that all vertices represent multicurves without null homologous submulticurves.\\

Suppose now that $\beta_i$ contains $k$ curves in the homotopy class $\left[b\right]$, where $k>1$. Assume also that $\beta_i$ and $m_2$ are in minimal position. The overlap function $f_i$ of $\beta_i$ and $m_2$ increases by one when crossing over a curve in $\left[b\right]$ from left to right. The variation (i.e. the maximum minus the minimum) of $f_i$ over all subsurfaces of $S$ adjacent to any curve of $\beta_i$ in $\left[b\right]$ is therefore equal to the variation of $f_i$ over all subsurfaces of $S$ adjacent to a fixed curve of $\beta_i$ in $\left[b\right]$, plus $k-1$. The possibility has not been ruled out that the existence of $k$ homotopic curves in $\beta_i$ might make it possible to construct $\beta_{i+1}$ with $\delta(\beta_{i},m_{2})=\delta(\beta_{i+1},m_{2})-3\chi(S)+k-2$. However, if $\beta_i$ has $k$ homotopic curves, it follows that there were at least $k-1$ surgeries performed at some stage in the construction of the path $m_{1}, \beta_{1}, \beta_{2},\ldots, \beta_i$ that did not give rise to null homologous submulticurves that could each be discarded to reduce the homological distance by one. Therefore, the decrease in homological distance from $m_2$ was overestimated by at least $k-1$ in previous steps, so the average over $i$ of $\delta(\beta_{i},m_{2})-\delta(\beta_{i+1},m_{2})$ is still bounded from above by $-3\chi(S)$.

\end{proof}




\begin{figure}
\centering
\includegraphics[width=13cm]{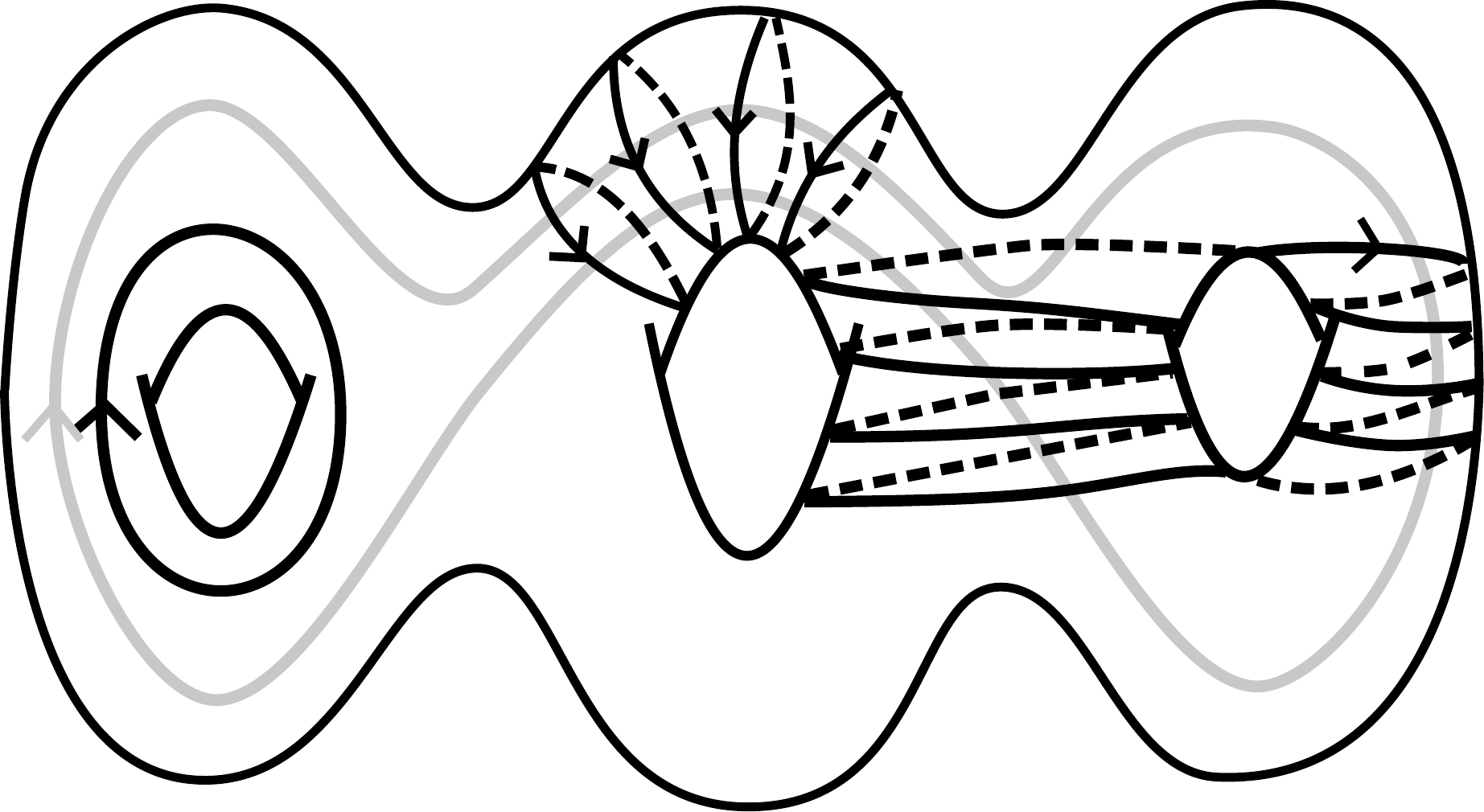}
\caption{A multicurve $m$ homologous to a simple curve (drawn in grey). The multicurve $m$ contains homotopic curves and no null homologous submulticurves.}
\label{chtwo7}
\end{figure}

\textbf{Remark.} The assumption that $\alpha$ is primitive is necessary in the proof of theorem \ref{quasigeodesic}. If $\alpha$ were not a primitive homology class, for example, if $\alpha$ is homologous to $nm_1$, $\delta(nm_{1},nm_{2})=n\delta(m_{1},m_{2})$, although the distance between $nm_1$ and $nm_2$ in $\mathcal{HC}(S, [nm_{1}])$ is equal to the distance between $m_1$ and $m_2$ in $\mathcal{HC}(S, [m_{1}])$.

\section{Minimal Genus Surfaces.}
\label{surfacesection}
This section gives a proof of theorem \ref{afterthought} and an algorithm for constructing minimal genus surfaces.\\

An aim of this paper is to use paths in $\mathcal{HC}(S,\alpha)$ to describe surfaces with boundary in $S\times \mathbb{R}$. Theorem 1.1 of \cite{Me2} states that every embedded, oriented, incompressible surface in $S\times I$ with boundary $m_{2}-m_{1}$ is constructed from a path in $\mathcal{HC}(S,\alpha)$. We want to use this theorem to describe minimal genus surfaces. In order to do this, it is necessary to establish whether or not a minimal genus surface is necessarily homotopic to an embedded surface. The answer to this is no; as shown in example \ref{reconnect}. However, it will be shown that there always exists an embedded surface with minimal genus. \\

The reason that the existence of an embedded, minimal genus surface is not immediately clear is due to possible intersections of the boundary of a surface with its interior. In this case there is no obvious surgery to remove such intersections without changing the boundary of the surface. It is implicit in the proof of theorem \ref{afterthought} that a minimal genus surface is homotopic to a surface whose boundary is disjoint from its interior.\\



\textbf{Euler Integrals.} Euler characteristic satisfies the properties of a measure, for example, for submanifolds $A$ and $B$ of $S\times 0$, $\chi(A\cup B)=\chi(A)+\chi(B)-\chi(A\cap B)$. Integration with respect to Euler characteristic, defined in \cite{Rota}, is a homomorphism from the ring of integer valued functions into $\mathbb{Z}$ such that
\begin{equation*}
\int 1_{A}d\chi =\chi(A)
\end{equation*}
Where $1_{A}$ is the function equal to 1 on the set $A$ and zero elsewhere. Integration with respect to Euler characteristic is defined here to relate $g_{H}$ to $\chi(H)$. Euler integrals have recently been used in a similar way for performing counting arguments on a complex that arises from studying sensor networks, \cite{Ghrist}.

It follows directly from the definition that for each $i$,
\begin{equation}
\label{Euler}
\int_{S\times \left\{0\right\}}|g_{T_{i}}|d\chi = \chi(T_{i})
\end{equation}
where the $\left\{T_i\right\}$ are defined in section \ref{surfaceconstruction}.\\

We now have all the necessary ingredients for proving theorem \ref{afterthought}.

\begin{proof}{Proof of theorem \ref{afterthought}}
Let $\mathcal{F}$ be the set of all oriented surfaces in $S\times I$ with boundary $m_{2}-m_{1}$, for homologous multicurves $m_1$ and $m_2$. Let $H$ be a surface in $\mathcal{F}$ with minimal genus.

Since the Euler characteristic of a cylinder is zero, it follows from lemma 6.1 of \cite{Me2} that
\begin{equation}
\label{equality}
\left|\chi(H)\right|=\sum_{i}\left|\chi(T_{i})\right|=\sum_{i}\int_{S\times \left\{0\right\}}\left|g_{T_{i}}\right|d\chi \geq \int_{S\times \left\{0\right\}}\left|g_{H}\right|d\chi
\end{equation}
Consider all subsurfaces $\pi(T_{i})\cap \pi(T_{j})$ of $S\times 0$ with nonzero Euler characteristic. Equality is achieved in the above equation iff all such intersecting subsurfaces are oriented in the same way.\\


Let $\gamma$ be the path constructed by the path construction algorithm. A surface $M$ is constructed from $\gamma$, such that $\chi(M)=\int_{S\times \left\{0\right\}}|g_{H}|d\chi$. This can be achieved by a specific choice of the component of $S\setminus (\gamma_{i+1}-\gamma_{i})$ the subsurface $T_i$ (defined in subsection \ref{surfaceconstruction}) should be homotopic to. 

By convention, the subsurface of $S\times \left\{0\right\}$ to the left of $\gamma_{i+1}-\gamma_i$ is oriented as a subsurface of $S\times \left\{0\right\}$, and the subsurface of $S\times \left\{0\right\}$ to the right of $\gamma_{i+1}-\gamma_i$ has the opposite orientation. Suppose the maximum of $g_{H}$ is equal to $k$. Then for $i\leq k$, choose $T_i$ to be to the left of $\gamma_{i+1}-\gamma_i$, and for all other $i$, choose $T_i$ to be to the right of $\gamma_{i+1}-\gamma_i$. It follows by construction that $g_{M}=g_{H}$.\\

Let $C_i$ denote the subsurface of $S\times \left\{0\right\}$ to the left of $\gamma_{i+1}-\gamma_i$. As discussed in the proof of lemma \label{gp}, the multicurves $\left\{\gamma_{i}\right\}$ obtained from the path construction algorithm are representatives of their homotopy classes such that $C_{1} \subset C_{2} \subset C_{3} \subset \ldots $ For all $i > k$ and $j\leq k$. It follows that $\pi(T_{i})$ is disjoint from $\pi(T_{j})$, and equality is achieved in equation \ref{equality}.\\

Since $M$ is constructed from a simple path in $\mathcal{HC}(S,\alpha)$, it is embedded.
\end{proof}

\textbf{Construction of minimal genus surfaces.} 

By equation \ref{equality}, constructing a minimal genus surface is reduced to the problem of finding the constant $c$ for which $\int_{S\times \left\{0\right\}}\left|g_{H}+c\right|d\chi$ is minimised. A necessary condition for $\int_{S\times \left\{0\right\}}\left|g_{H}+c\right|d\chi$ to be minimised is that $-\delta(m_{1},m_{2}) \leq g_{H}+c \leq \delta(m_{1},m_{2})$. This leaves a finite number of choices for c.\\

Actually, there is a uniform bound on the number of choices for $c$. Due to the fact that $C_{1} \subset C_{2} \subset C_{3} \subset \ldots $, Euler characteristic arguments give a uniform bound on the number of indices $i$ such that $\chi(C_{i})\neq \chi(C_{i+1})$. Let $k$ be the smallest value of $i$ such that $\left|\chi(S\setminus C_{i+1})\right| <\left|\chi(C_{i+1})\right|$. A minimal genus surface $M$ is obtained by choosing $T_i$ to be to the left of $\gamma_{i+1}-\gamma_i$ for $i\leq k$, and for all other $i$, $T_i$ is to the right of $\gamma_{i+1}-\gamma_i$. It is clear that $M$ has minimal genus, since equality is achieved in equation \ref{equality}, and the choice of $k$ ensures that $\int_{S\times \left\{0\right\}}|g_{M}|d\chi$ is minimised.\\

\begin{ex}[Example of a minimal genus surface that is not embedded]
\label{reconnect}
Let $m_1$ and $m_2$ be homologous curves as shown in figure \ref{spiralwithears}. Let $H$ be the embedded, minimal genus surface with boundary $m_{2}-m_{1}$ constructed as in the previous paragraph. The values of $g_{H}$ are shown in the figure. Let $c$ be a simple oriented curve contained in a subsurface of $S\times \left\{0\right\}$ on which $g_{H}\geq 2$. 
Suppose $c_1$ and $c_2$ are two curves in $H$ that are both homotopic to the curve $c$ in $S\times \left\{0\right\}$. Cutting $H$ along $c_1$ and $c_2$ gives a surface with boundary $m_{2}-m_{1}+c_{1}-c_{1}+c_{2}-c_{2}$. Construct a surface $H^{'}$ by gluing the boundary curve $c_1$ to $-c_2$ and then gluing $c_{2}$ to $-c_1$. Clearly, $H^{'}$ has the same boundary, Euler characteristic and pre-image function as $H$, but can not be embedded.
\begin{figure}
\begin{center}
\def\svgwidth{12cm}
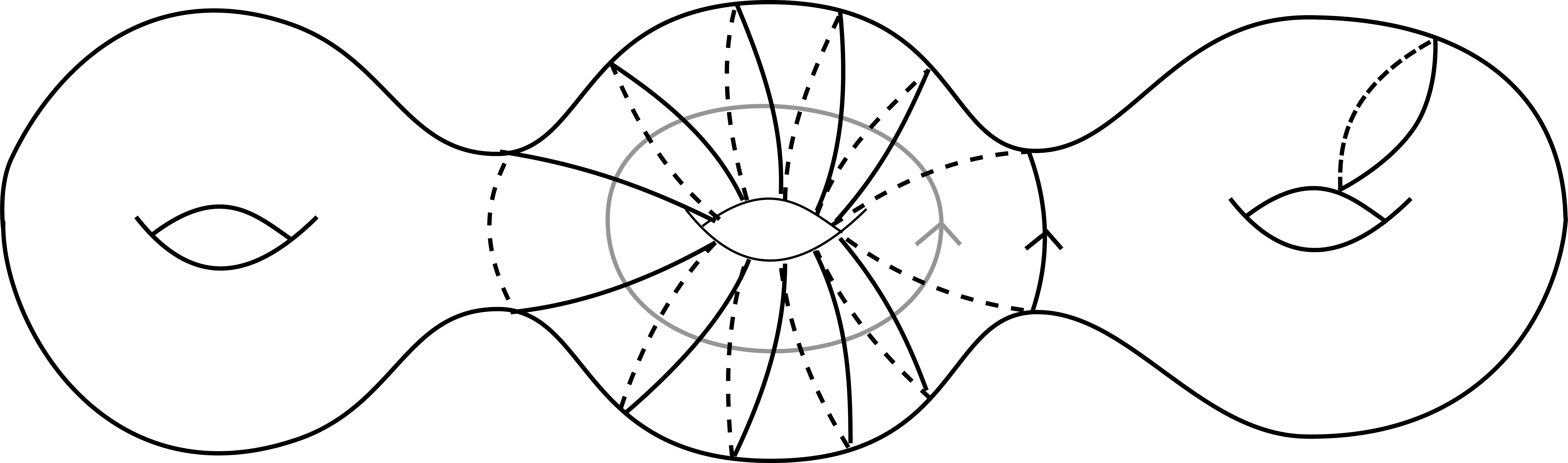
\label{spiralwithears}
\end{center}
\end{figure} 

\end{ex}

\section{Quasi-flats and Distance Bounds}
\label{Examples}
This section gives a few simple examples to illustrate key geometric properties of $\mathcal{HC}(S,\alpha)$.\\

In theorem \ref{j}, distances in $\mathcal{HC}(S,\alpha)$ were shown to be related to the homological distance $\delta$. The next question is, how does distance relate to intersection number? At each step of the path construction algorithm, the intersection number with $m_2$ is decreased. Recall that the arcs of $m_{2}\cap (S\setminus m_{1})$ on $\partial S_{max}(m_{1},m_{2})$ were denoted $a_{1},\ldots, a_n$. Let $k_{a_i}$ be the number of arcs of $m_{2}\cap (S\setminus m_{1})$ in the same homotopy class as $a_i$ for $1\leq i \leq n$. Then the intersection number of $\gamma_1$ with $m_2$ is at least $2\sum_i k_{a_i}$ less than the intersection number of $m_1$ with $m_2$. \\

It is well known that the distance between two curves $c_1$ and $c_2$ in the complex of curves is bounded from above by $a\log(i(c_{1},c_{2}))+b$, for some constants $a$ and $b$. In figure \ref{chthreeone1} an example is given that demonstrates that the distance between two curves in $\mathcal{HC}(S,\alpha)$ can be as much as $\frac{i(c_{1},c_{2})}{2}+1$. \\

\begin{proof}[Proof of theorem \ref{thmtwo}]
\label{toruswithears}
Let $c_{1}$ and $c_2$ be the curves shown in figure \ref{chthreeone1}. The curve $c_{2}$ is obtained by Dehn twisting $c_{1}$ $n$ times around one curve, $t_1$, in a bounding pair, and $-n$ times around the other curve, $t_2$, in the bounding pair. In figure \ref{chthreeone1}, $n$ is five. A simple calculation shows that $\delta(c_{1},c_{2})$ is equal to $\frac{i(c_{1},c_{2})}{2}+1$. 

To see why the distance between $c_1$ and $c_2$ can't be less than $\frac{i(c_{1},c_{2})}{2}+1$, observe that any multicurve in $\alpha$ has nonzero algebraic intersection number with each of $t_1$ and $t_2$. Suppose $c_{1}, \beta_{1},\ldots, c_{2}$ is a path connecting $c_1$ and $c_2$ in $\mathcal{HC}(S,\alpha)$. Informally, it follows that $i(\beta_{i},c_{2}) \geq i(c_{1},c_{2})-2i$, because it is not possible to unwind more than one pair of twists at each step. To be more precise, in \cite{MasurandMinskyII}, distance between $c_1$ and $c_2$ in the subsurface projection to an annulus with core curve $c$ was defined. The distance between $c_1$ and $c_2$ in the subsurface projection to an annulus with core curve $t_1$ depends linearly on the intersection number of the lifts of $c_1$ and $c_2$ to the covering space consisting of an annulus with core curve $t_1$. Unlike in the complex of curves, a path in $\mathcal{HC}(S,\alpha)$ has to pass through the subsurfaces consisting of annuli whose core curves have nonzero algebraic intersection number with $\alpha$.


\begin{figure}
\centering
\includegraphics[width=13cm]{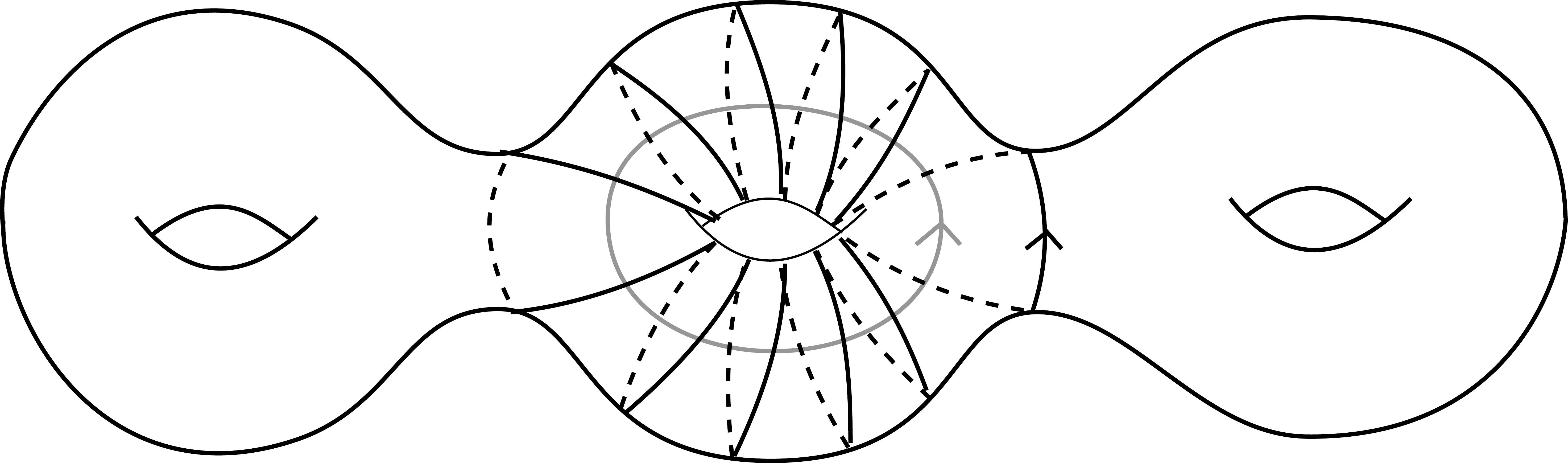}
\caption{Example demonstrating that the best possible upper bound on the distance between $c_{1}$ and $c_{2}$ in $\mathcal{HC}(S, \alpha)$ is given by $\frac{i(c_{2},c_{1})}{2}+1$.} 
\label{chthreeone1}
\end{figure}
 
\end{proof}

\begin{figure}
\begin{center}
\def\svgwidth{13cm}
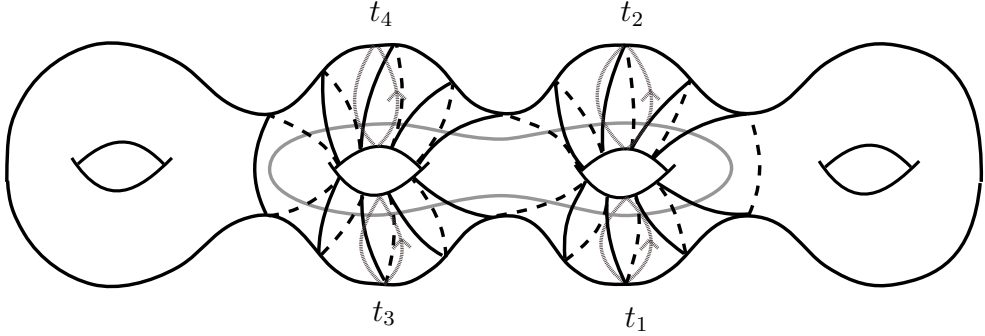
\end{center}
\label{doubletoruswithears}
\caption{How to construct a quasi-flat.}
\end{figure}

The example shown in figure 13 is generalised to construct families of examples to show that $\mathcal{HC}(S,\alpha)$ is not $\delta$-hyperbolic for $g>3$. 

\begin{proof}[Proof of theorem \ref{thmone}]
For $g>3$ there exist two pairs of bounding pairs $(t_{1}, -t_{2})$ and $(t_{3},-t_{4})$; each of the $t_i$ representing distinct homotopy classes. Suppose $v_{1}$ is a multicurve with nonzero algebraic intersection number with each of $t_{1}, t_{2}, t_{3}$ and $t_{4}$, as in figure 13. Let $v_2$ be the multicurve $v_1$ Dehn twisted around $(t_{1}, t_{2})$ $n$ times, and $v_{3}$ be the multicurve $v_{1}$ Dehn twisted around $(t_{3},t_{4})$ $n$ times. The vertices of a geodesic triangle in $\mathcal{HC}(S,\alpha)$ are represented by the symbols $v_1$, $v_2$ and $v_3$. 
Since the distance between two vertices on the boundary of the triangle is equal to the number of Dehn twists around $(t_{1}, t_{2})$ and $(t_{3},t_{4})$ necessary to get from one vertex to the other, for $n$ even, the midpoints of the sides of the geodesic triangle are each a distance $\frac{n}{2}$ from the other two sides of the triangle. The number of twists, $n$, can therefore be chosen large enough so that this triangle is not $\delta$-thin. 
\end{proof}

\textbf{Tightness.} Curve complexes are in general locally infinite, so there can be infinitely many geodesic paths connecting two vertices. However, most of these geodesics do not seem to provide any additional structural information. In order to be able to prove finiteness results in a locally infinite complex, the concept of ``tightness'' was introduced in \cite{MasurandMinskyII}. The definition given here is from \cite{Bowditch3}. A path $\gamma_{1}, \gamma_{2},\ldots,\gamma_n$ in $\mathcal{HC}(S,[\gamma_{1}])$ is \textit{tight} at some index $i\neq 1,n$ if every closed curve that intersects $\gamma_i$ also intersects $\gamma_{i-1}\cup \gamma_{i+1}$. The path is tight if it is tight for all $i\neq 1,n$. It follows from lemma \ref{gp} that all paths constructed by the path construction algorithm are tight.

The example in figure 13 also shows that, unlike in the complex of curves, in $\mathcal{HC}(S,\alpha)$ there does not always exist a tight geodesic connecting any two vertices. A geodesic $c_{1}, \gamma_{1}, \gamma_{2}, \ldots, c_{2}$ is constructed such that for each $i$, $\gamma_{i+1}$ is obtained from $\gamma_i$ by performing Dehn twists around $t_1$, $t_2$, $t_3$ and $t_4$. It is not hard to check that this is only possible if $\gamma_1$ is obtained from $c_{1}$ by performing a surgery that cuts $c_{1}$ into two curves; one that intersects $t_1$ and $t_2$, and another one that intersects $t_3$ and $t_4$, see for example the middle diagram of figure \ref{doubletoruswithears2}. All curves contained in the one dimensional cell complex $c_{1}\cup c_{2}$ are either null homologous, $c_{1}$, $c_{2}$, $t_1$, $t_2$, $t_3$, $t_4$ or they intersect all of $t_1$, $t_2$, $t_3$ and $t_4$. It follows that a geodesic connecting $c_{1}$ and $c_{2}$ can not be tight.


\section{Appendix - Distances in the Cyclic Cycle Complex.}
\label{Bike}
In \cite{Hatcher}, a closely related complex, the cyclic cycle complex, was defined. In this appendix, it is shown that the path construction algorithms from section \ref{algorithm} can be modified slightly to construct geodesics in this complex.\\

A multicurve $m$ is said to be \textit{reduced} if it does not contain a submulticurve that bounds a complementary region of $m$ in $S$ (using either orientation of the region). The \textit{Cyclic Cycle Complex} $\mathcal{CC}(S)$ from \cite{Hatcher} is the simplicial complex whose vertices are the homotopy classes of oriented, reduced multicurves. A set of $k+1$ vertices spans a simplex in $\mathcal{CC}(S)$ if these vertices are represented by disjoint multicurves $m_0$, $m_1$, $m_{2}$,\ldots,$m_k$ that cut $S$ into $k+1$ embedded subsurfaces $E_0$, $E_1$,\ldots,$E_{k}$ such that the oriented boundary of $E_i$ is $m_{i+1}-m_{i}$. In particular, all edges are by definition simple. As a consequence, paths in the cyclic cycle complex correspond to embedded surfaces in $S\times I$, as opposed to merely immersed.\\

It follows that each connected component of $\mathcal{CC}(S)$ represents multicurves in a fixed nontrivial homology class. Every connected component of $\mathcal{CC}(S)$ can therefore be embedded in $\mathcal{HC}(S,\alpha)$ for appropriate $\alpha$.\\

\begin{thm}
Let $m_1$ and $m_2$ be multicurves representing vertices in the same connected component of $\mathcal{CC}(S)$. The distance between $m_1$ and $m_2$ in $\mathcal{CC}(S)$ is equal to $\delta(m_{1},m_{2})$, and geodesic paths can be explicitly constructed.
\label{modified}
\end{thm}
\begin{proof}
The path construction algorithm can be easily modified to construct paths in $\mathcal{CC}(S)$. A vertex in $\mathcal{HC}(S, \alpha)$ might not correspond to a vertex in $\mathcal{CC}(S)$, because the multicurves representing vertices in $\mathcal{CC}(S)$ are not allowed to contain just any null homologous submulticurve. Let $\gamma_{1}, \gamma_{2},\ldots,\gamma_k$ be the simple path in $\mathcal{HC}(S,\alpha)$ constructed by the path construction algorithm. Suppose also that $\gamma_1$ and $\gamma_k$ represent reduced multicurves. A path $\gamma_{1}^{'}, \gamma_{2}^{'}$,\ldots,$\gamma_{k}^{'}$ in $\mathcal{CC}(S)$ can be constructed as follows: Whenever the vertex $\gamma_{i}$ corresponds to a multicurve $\tilde{\gamma}_{i}$ containing a null homologous submulticurve $n$ that bounds a complementary region of $\tilde{\gamma}_{i}$ in $S$, (i.e. $\tilde{\gamma}_{i}$ is not reduced) let  $\tilde{\gamma}_{i}^{'}$ be the multicurve $\tilde{\gamma}_{i}-n$. If $\tilde{\gamma}_{i}$ is reduced, let $\tilde{\gamma}_{i}=\tilde{\gamma}_{i}^{'}$. The symbol $\gamma_{1}^{'}$ denotes the vertex of $\mathcal{CC}(S)$ corresponding to the multicurve $\gamma_{i}^{'}$. It remains to show that $\gamma_{i+1}^{'}$ and $\tilde{\gamma}_{i}^{'}$ are connected by an edge in $\mathcal{CC}(S)$ for all $1\leq i \leq k-1$.

Let $j$ be the smallest integer such that $\tilde{\gamma}_j$ is not reduced, and let $n$ be the union of null homologous submulticurves of $\tilde{\gamma}_j$ that bound complementary regions of $\tilde{\gamma}_j$ in $S$. To show that $\tilde{\gamma}_{j}^{'}$ and $\tilde{\gamma}_{j-1}^{'}$ are connected by an edge in $\mathcal{CC}(S)$ consists of showing that there is a subsurface $N$ of $S$ (with either orientation) with $\partial N=n$ and such that $N$ is disjoint from $S_{max}(\tilde{\gamma}_{j-1}^{'},\tilde{\gamma}_{j}^{'})$.


By construction $S_{max}(\tilde{\gamma}_{j-1}^{'}, \tilde{\gamma}_{j}^{'})$ and $\tilde{\gamma}_{j}^{'}-\tilde{\gamma}_{j-1}^{'}$ is to the left of $n$ (If $S_{min}$ had been used in place of $S_{max}$ in the path construction algorithm, $\tilde{\gamma}_{j}^{'}-\tilde{\gamma}_{j-1}^{'}$ would have to be to the right of $n$). Also, since $\tilde{\gamma}_{j-1}^{'}$ is reduced, it follows from the arguments given in the proof of theorem \ref{quasigeodesic} that $\tilde{\gamma}_{j}^{'}$ can not contain a second null homologous submulticurve $n_2$ that lies between $n$ and the other curves in $\tilde{\gamma}_{j}^{'}-\tilde{\gamma}_{j-1}^{'}$. The null homologous multicurve $n$ can therefore be capped off from the right by a subsurface disjoint from $S_{max}(\tilde{\gamma}_{j-1}^{'},\tilde{\gamma}_{j}^{'})$.

Suppose the null homologous multicurve $n$ from the previous paragraph is a submulticurve of $\tilde{\gamma}_i$ but not $\tilde{\gamma}_{i}^{'}$ for some $j<i$, and the surgery performed on $\tilde{\gamma}_{i}$ to obtain $\tilde{\gamma}_{i+1}$ alters the submulticurve $n$. Then the subsurface of $S$ bounded by $\tilde{\gamma}_{i+1}^{'}-\tilde{\gamma}_{i}^{'}$ contains the subsurface $N$ of $S$. A symmetric argument, in which $m_1$ and $m_2$ are exchanged shows that the subsurface of $S$ bounded by $\tilde{\gamma}_{i+1}^{'}-\tilde{\gamma}_{i}^{'}$ is still an embedded subsurface of $S$.

The claim then follows by induction.

That the path in $\mathcal{CC}(S)$ so constructed is a geodesic follows from theorem \ref{j}.
\end{proof}




\label{Hatcheralgo}

\bibliographystyle{plain}
\nocite{*}
\bibliography{bibfile}
\end{document}